\documentclass[11pt]{article}
\usepackage{amssymb}
\usepackage{amscd}
\usepackage{times}
\usepackage{graphics}
\usepackage{amscd}
\usepackage{amsmath}
\usepackage[all]{xypic}

\begin{document}
\newtheorem{theorem}{Theorem}[section]
\newtheorem{lemma}[theorem]{Lemma}
\newtheorem{corollary}[theorem]{Corollary}
\newtheorem{definition}[theorem]{Definition}
\newtheorem{proposition}[theorem]{Proposition}
\newtheorem{defprop}[theorem]{Definition-Proposition}
\newtheorem{example}[theorem]{Example}
\newtheorem{remark}[theorem]{Remark}

\newcommand{\pf}{\noindent {\bf Proof:} }
\catcode`\@=11
\@addtoreset{equation}{section}
\catcode`\@=12

\renewcommand{\theequation}{\arabic{section}.\arabic{equation}}
\def\lu{\rightharpoonup}
\def\sqr#1#2{{\vcenter{\vbox{\hrule height.#2pt\hbox{\vrule width.#2pt
  height#1pt \kern#1pt \vrule width.#2pt}\hrule height.#2pt}}}}
\def\square{\mathchoice\sqr64\sqr64\sqr{2.1}3\sqr{1.5}3}
\def\argh{\rightharpoonup}
\def\m#1{m_{(#1)}}
\def\g#1{g_{(#1)}}
\def\h#1{h_{(#1)}}
\def\k#1{k_{(#1)}}
\def\M{{\cal M}}
\def\Z{{\mathbb Z}}
\def\F{{\cal F}}
\def\mm#1{m_{(#1*)}}
\def\a#1{a_{(#1)}}
\def\aa#1{a_{(#1*)}}
\def\b#1{b_{(#1)}}
\def\bb#1{b_{(#1*)}}
\def\l#1{l_{(#1)}}
\def\n#1{n_{(#1)}}
\def\nn#1{n_{(#1*)}}
\def\totimes{\tilde{\otimes}}
\def\c#1{c_{(#1)}}
\def\cc#1{c_{(#1*)}}
\def\id{{\rm id}}
\def\ydh{{\cal YD}(H)}
\def\yddh{{\cal YD}(H^*)}
\def\lefth{_{H}{\cal M}}
\def\lefthd{_{H^*}{\cal M}}
\def\righth{{\cal M}^H}
\def\righthd{{\cal M}^{H^*}}
\def\C{{\cal C}}
\def\D{{\cal D}}
\def\R{{\cal R}}
\def\rsb{r_{s,\beta}}
\def\rsg{r_{s,\gamma}}
\def\km{k^{*}}
\def\pf{{\bf Proof:} }
\def\U{{\cal U}}
\def\T{{\cal T}}
\def\H{{\cal H}}
\def\Zen{Z_L(E(n))}
\def\sym{Sym_{M,n,r}(k)}
\def\ru{\rightharpoonup}

\title{On the subgroup structure of the full Brauer group of Sweedler Hopf algebra}
\author{{\small
\begin{tabular}{p{6.5cm}p{0.0cm}p{7cm}}
{\Large Giovanna Carnovale} &  &  \hspace{-0.7cm} {\Large Juan Cuadra} \vspace{4pt}\\
Dipartimento di Matematica Pura  &  & \hspace{-0.7cm} Universidad de Almer\'{\i}a \\
ed Applicata& &  \hspace{-0.7cm} Dpto. \'{A}lgebra y An\'{a}lisis Matem\'atico \\
via Trieste 63 & & \hspace{-0.7cm} E-04120 Almer\'{\i}a, Spain \\
I-35121 Padua, Italy &  & \hspace{-0.7cm} jcdiaz@ual.es \\
carnoval@math.unipd.it & &
\end{tabular}
}}

\date{}
\maketitle

\begin{abstract}
We introduce a family of three parameters $2$-dimensional algebras representing elements in the Brauer group
$BQ(k,H_4)$ of Sweedler Hopf algebra $H_4$ over a field $k$. They allow us to describe the mutual intersection of the subgroups arising from a quasitriangular or coquasitriangular structure. We also define a new subgroup of $BQ(k,H_4)$ and construct an exact sequence relating it to the Brauer group of Nichols $8$-dimensional Hopf algebra with respect to the  quasitriangular structure attached to the $2\times 2$-matrix with $1$ in the $(1,2)$-entry and zero elsewhere.
\end{abstract}


\section*{Introduction}

The Brauer group of a Hopf algebra is an extremely complicated
invariant that reflects many aspects of the Hopf algebra: its
automorphisms group, its Hopf-Galois theory, its second lazy
cohomology group, (co)quasitriangularity, etc.  It is very difficult
to describe all its elements and to find their multiplication rules.
For the most studied case, that of a commutative and cocommutative
Hopf algebra, these are the results known so far: the first explicit
computation was done by Long in \cite{L1} for the group algebra
$k\Z_n,$ where $n$ is square-free and $k$ algebraically closed with
$char(k) \nmid n$; DeMeyer and Ford \cite{DF} computed it for $k\Z_2$
with $k$ a commutative ring containing $2^{-1}$. Their result was
extended by Beattie and Caenepeel in \cite{BC} for $k\Z_n,$ where $n$
is a power of an odd prime number and some mild assumptions on $k$.
In \cite{C1} Caenepeel achieved to compute the multiplication rules
for a subgroup, the so-called split part, of the Brauer group for a
faithfully projective commutative and cocommutative Hopf algebra $H$
over any commutative ring $k$. These results were improved in
\cite{C2} and allowed him to compute the Brauer group of Tate-Oort
algebras of prime rank. For a unified exposition of these results
the profuse monograph \cite{C3} is recommended.
\smallskip

Since the Brauer group was defined for any Hopf algebra with bijective
antipode (\cite{CVZ}, \cite{CVZ2}), it was a main goal to compute it
for the smallest noncommutative noncocommutative Hopf algebra:
Sweedler's four dimensional Hopf algebra $H_4$, which is generated
over the field $k$ ($char(k) \neq 2$) by the group-like $g$, the
$(g,1)$-primitive element $h$ and relations $g^2=1, h^2=0, gh=-hg$. A first step was the calculation in \cite{VOZ3} of the subgroup $BM(k,H_4,R_0)$ induced by the quasitriangular structure $R_0=2^{-1}(1 \otimes 1+g\otimes 1+1 \otimes g-g \otimes g).$ It was shown to be isomorphic to the direct product of $(k,+)$, the additive group of $k$, and $BW(k)$, the Brauer-Wall group of $k$. It was later proved in \cite{gio} that the subgroups $BM(k,H_4,R_t)$ and $BC(k,H_4,r_s)$  arising from all the quasitriangular structures $R_t$ and the coquasitriangular structures $r_s$ of $H_4$ respectively, with $s,t \in k$, are all isomorphic. \smallskip

In this paper we introduce a family of three parameters $2$-dimensional algebras $C(a;t,s)$, for $a,t,s \in k,$ that represent elements in $BQ(k,H_4)$. They will allow us to shed a ray of light on the subgroup structure of $BQ(k,H_4)$ and will provide some evidences about the difficulty of the computation of this group. The algebra $C(a;t,s)$ is generated by $x$ with relation $x^2=a$ and has a $H_4$-Yetter-Drinfeld module algebra structure with action and coaction:
$$g\cdot x=-x,\quad h\cdot x=t, \quad \rho(x)=x\otimes g+s\otimes h.$$
We list the main properties of these algebras in Section 2 (Lemma \ref{properties}) and we show that $C(a;t,s)$ is $H_4$-Azumaya if and only if $2a \neq st$. When $s=lt$ they represent elements in $BM(k,H_4,R_l)$ and this subgroup is indeed generated by the classes of $C(a;1,t)$ with $2a \neq t$ together with $BW(k)$, Proposition \ref{BMt}. The same statement holds true for
$BC(k,H_4,r_l)$ when $t=sl$ replacing $C(a;1,t)$ by $C(a;s,1)$, Proposition \ref{BCs}. \smallskip

Using the description of $BM(k, H_4, R_t)$ and $BC(k, H_4, r_s)$ in terms of these algebras, Section 3 is devoted to analyze the intersection of these subgroups inside $BQ(k, H_4)$. Let $i_t$ and $\iota_s$ denote the inclusion map of the former and the latter respectively. It is known that $BW(k)$ is contained in any of the above subgroups. Theorem \ref{intersection} states that:
\begin{enumerate}
\itemsep -1pt
\item[(1)] $Im(i_t)\cap Im(\iota_s)\neq BW(k)$ iff $ts=1$. If this is the case, $Im(i_t)=Im(\iota_s)$;
\item[(2)] $Im(i_t)\cap Im(i_s)\neq BW(k)$ if and only if $t=s$;
\item[(3)] $Im(\iota_t)\cap Im(\iota_s)\neq BW(k)$ if and only if $t=s$.
\end{enumerate}
A remarkable property of our algebras is that they represent the same class in $BQ(k, H_4)$ if and only if they are isomorphic, Corollary \ref{equal}. \smallskip

A morphism from the automorphism group of $H_4$ to $BQ(k,H_4)$ was constructed in \cite{VOZ2}, allowing to consider
$k^{\cdot 2}$ as a subgroup of $BQ(k,H_4)$. In Section 4 we show that the subgroup $BM(k,H_4, R_l)$ is conjugated to $BM(k,H_4, R_{l\alpha^2})$ inside $BQ(k, H_4)$, for $\alpha \in k^\cdot$, by a suitable representative of $k^{\cdot 2}$, Lemma
\ref{conj-autom}. \smallskip

Any $H_4$-Azumaya algebra possesses two natural ${\mathbb Z}_2$-gradings: one stemming from the action of $g$
and one from the coaction (after projection) of $g$. In Section 6 we introduce the subgroup
$BQ_{grad}(k, H_4)$ consisting of those classes  of $BQ(k,H_4)$  that can be represented by $H_4$-Azumaya algebras for which the two ${\mathbb Z}_2$-gradings coincide. On the other hand, the Drinfeld double of $H_4$ admits a Hopf algebra map $T$ onto Nichols $8$-dimensional Hopf algebra $E(2)$. This map is quasitriangular as $E(2)$ is equipped with the quasitriangular structure $R_N$ corresponding to the $2\times 2$-matrix $N$ with $1$ in the $(1,2)$-entry and zero elsewhere, see (\ref{RN}). If we consider the associated Brauer group $BM(k,E(2),R_N)$, then Theorem \ref{sequence} claims that $T$ induces a group homomorphism $T^*$ fitting in the following exact sequence
$$\begin{array}{l}
\CD 1\longrightarrow{\mathbb Z}_2@>>>BM(k,
E(2),R_N)@>{T^*}>>BQ_{grad}(k,H_4)\longrightarrow 1.
\endCD
\end{array}$$
So in order to compute $BQ(k,H_4)$ one should first understand $BM(k,E(2),R_N)$. This new problem cannot be attacked with the available techniques for computations of groups of type BM, \cite{VOZ3}, \cite{hnu}, \cite{enne}. Those computations were achieved by finding suitable invariants for a class by means of a Skolem-Noether-like theory. In the Appendix we underline some obstacles to the application of these techniques to the computation of $BM(k,E(2),R_N)$: the set of elements represented by algebras for which the action of one of the standard nilpotent generators of $E(2)$ is inner coincides with the set of classes represented by ${\mathbb Z}_2$-graded central simple algebras and this is not a subgroup of $BM(k,E(2),R_N),$ Theorems \ref{conditions}, \ref{not-subgroup}. Moreover, $BM(k,E(2),R_N)$ seems to be much more complex than the groups of type BM treated until now since, according to Proposition \ref{lambdamu}, each group $BM(k,H_4,R_t)$ may be viewed as a subgroup of it.

\section{Preliminaries}\label{preliminaries}

In this paper $k$ is a field, $H$ will denote a Hopf algebra over $k$ with bijective antipode $S$, coproduct $\Delta$ and counit $\varepsilon$. Tensor products $\otimes$ will be over $k$ and, for vector spaces $V$ and
$W$, the usual flip map is denoted by $\tau: V\otimes W\to W\otimes V$.  We shall adopt the Sweedler-like notations $\Delta(h)=\h1\otimes \h2$ and $\rho(m)=\m0\otimes \m1$ for coproducts and right comodule structures respectively. For $H$ coquasitriangular (resp. quasitriangular), the set of all coquasitriangular (resp. quasitriangular) structures will be denoted by $\cal U$ (resp. $\cal T$). \par \bigskip

{\it Yetter-Drinfeld modules.} Let us recall that if $A$ is a left $H$-module with action $\cdot$ and a right $H$-comodule with coaction $\rho$ the two structures combine to a left module structure for the Drinfeld double $D(H)=H^{*, cop}\bowtie H$ of $H$ (cfr. \cite{M}) if and only if they satisfy the so-called Yetter-Drinfeld compatibility condition:
\begin{equation}\label{YD}
\rho(l\cdot b)=\l2\cdot\b0\otimes\l3\b1 S^{-1}(\l1), \quad \forall l\in H, b\in A.
\end{equation}
Modules satisfying this condition are usually called Yetter-Drinfeld modules. If $A$ is a left $H$-module algebra and a right $H^{op}$-comodule algebra satisfying (\ref{YD}) we shall call it a Yetter-Drinfeld $H$-module algebra. \medskip

{\it The Brauer group} (see \cite{CVZ}, \cite{CVZ2}). Suppose that $A$ is a Yetter-Drinfeld $H$-module algebra. The $H$-opposite algebra of $A$, denoted by $\overline{A}$, is the underlying vector space of $A$ endowed with product $a\circ c=\c0(\c1\cdot a)$ for every $a,c\in A$. The same action and coaction of $H$ on $A$ turn $\overline{A}$ into a Yetter-Drinfeld $H$-module algebra. Given two Yetter-Drinfeld $H$-module algebras $A$ and $B$ we can construct a new Yetter-Drinfeld module  $A\# B$ whose underlying vector space is $A\otimes B$, with action $h\cdot (a \otimes b)=\h1 \cdot a \otimes \h2 \cdot b$ and with coaction $a\otimes b\mapsto \a0\b0\otimes\b1\a1$. This object becomes a Yetter-Drinfeld module algebra if we provide it with the multiplication
$$(a\# b)(c\# d)=a\c0\#(\c1\cdot b) d.$$

For every finite dimensional Yetter-Drinfeld module $M$ the algebras ${\rm End}(M)$ and ${\rm End}(M)^{op}$ can be naturally provided of a Yetter-Drinfeld module algebra structure through (\ref{endm}) and (\ref{endmop}) below respectively:
\begin{equation}\label{endm}
\begin{array}{l}
(h \cdot f)(m)=h_{(1)} \cdot f (S(h_{(2)}) \cdot m), \vspace{2pt} \\
\rho(f)(m) = f(m_{(0)})_{(0)} \otimes S^{-1}(m_{(1)})f(m_{(0)})_{(1)},
\end{array}
\end{equation}
\begin{equation}\label{endmop}
\begin{array}{l}
(h \cdot f)(m)=h_{(2)} \cdot f (S^{-1}(h_{(1)}) \cdot m), \vspace{2pt} \\
\rho(f)(m) = f(m_{(0)})_{(0)} \otimes f(m_{(0)})_{(1)}S(m_{(1)}),
\end{array}
\end{equation}
where $h \in H, f \in End(M), m \in M.$ A finite dimensional Yetter-Drinfeld module algebra $A$ is called {\it $H$-Azumaya} if the following module algebra maps are isomorphisms:
\begin{equation}\label{f-and-g}
\begin{array}{ll}
F\colon A\#{\overline{A}} \rightarrow {\rm End}(A), &  F(a\# b)(c)=a\c0(\c1\cdot b), \vspace{2pt} \\
G\colon \overline{A}\#{{A}} \rightarrow {\rm End}(A)^{op}, & G(a\# b)(c)=\a0(\a1\cdot c)b.
\end{array}
\end{equation}
The algebras ${\rm End}(M)$ and ${\rm End}(M)^{op}$, for a finite dimensional Yetter-Drinfeld module $M$, provided with the preceding structures are $H$-Azumaya. \smallskip

The following relation $\sim$ established on the set of isomorphism classes of $H$-Azumaya algebras is an equivalence relation:
{\it $A\sim B$ if there exist finite dimensional Yetter-Drinfeld modules $M$ and $N$ such that $A\#{\rm End}(M)\cong B\# {\rm End}(N)$ as Yetter-Drinfeld module algebras}. The set of equivalence classes of $H$-Azumaya algebras, denoted by $BQ(k,H)$, is a group with product $[A][B]=[A\# B]$, inverse element $[\overline{A}]$ and identity element $[End(M)]$ for finite dimensional Yetter-Drinfeld modules $M$. This group is called the {\it full Brauer group of $H$}. The adjective full is used to distinguish it from the subgroups presented next, that receive the same name in the literature. \smallskip

Given a left $H$-module algebra $A$ with action $\cdot$ and a quasitriangular structure $R=R^{(1)}\otimes R^{(2)}$ on $H$, a right $H^{op}$-comodule algebra structure $\rho$ on $A$ is determined by
$$\rho(a)=(R^{(2)}\cdot a)\otimes R^{(1)}, \quad \forall a\in A.$$ We
will call this coaction the coaction induced by $\cdot$ and $R$. It is
well-known that $(A,\cdot,\rho)$ satisfies the Yetter-Drinfeld
condition. This allows the definition of the subgroup $BM(k,H,R)$ of
$BQ(k, H)$ whose elements are equivalence classes of $H$-Azumaya
algebras with coaction induced by $R$ (\cite[\S  1.5]{CVZ2}). To
underline that a representative $A$ of a given class in $BQ(k,H)$
represents a class in $BM(k,H,R)$ we shall say that {\it $A$ is an
  $(H,R)$-Azumaya algebra}. The inclusion map will be denoted by  $i\colon BM(k,H, R)\to BQ(k,H)$. For $H$ finite dimensional $BQ(k,H)=BM(k, D(H),{\cal R}),$ where ${\cal R}$ is the natural quasitriangular structure on the Drinfeld double $D(H)$.

Dually, given a right $H^{op}$-comodule algebra $A$ with coaction $\varrho$ and a coquasitriangular structure $r$ on $H$, a $H$-module algebra structure $\cdot$ on $A$ is determined by $$h\cdot a=\a0 r(h\otimes\a1), \quad \forall a\in A, h\in H,$$ and $(A,\cdot,\varrho)$ becomes a Yetter-Drinfeld module algebra. We will call this action the action induced by $\chi$ and $r$. The subset $BC(k,H,r)$ of $BQ(k, H)$ consisting of those classes admitting a representative whose action is induced by $r$ is a subgroup (\cite[\S  1.5]{CVZ2}). To stress that a representative $A$ of a class in $BQ(k,H)$ represents a class in $BC(k,H,r)$ we shall say that {\it $A$ is an $(H,r)$-Azumaya algebra}. The inclusion of $BC(k,H,r)$ in $BQ(k,H)$ will be denoted by  $\iota\colon BC(k,H, r)\to BQ(k,H)$. \bigskip

{\it On Sweedler Hopf algebra.} In the sequel we will assume that
$char(k) \neq 2.$ Let $H_4$ be Sweedler Hopf algebra, that is, the
Hopf algebra over  $k$ generated by a grouplike element $g$ and an
element $h$ with relations, coproduct and antipode:
$$g^2=1,\quad h^2=gh+hg=0, \quad \Delta(h)=1\otimes h+h\otimes g,\quad S(g)=g,\quad S(h)=gh.$$
The Hopf algebra $H_4$ has a family of  quasitriangular (indeed triangular) structures. They were classified in \cite{Rad} and are given by:
$$R_t=\frac{1}{2}(1\otimes 1+1\otimes g+g\otimes 1-g\otimes
g)+\frac{t}{2}(h\otimes h+h\otimes gh+gh\otimes gh-gh\otimes h),$$
where $t \in k$. It is well-known that $H_4$ is self-dual so that $H_4$ is also cotriangular. Let $\{1^*, g^*, h^*, (gh)^*\}$ be the basis of $H^*_4$ dual to $\{1,g,h,gh\}$. We will often
make use of the Hopf algebra isomorphism
$$\begin{array}{rl}
\phi\colon H_4&\to H_4^*\\
1&\mapsto 1^*+g^*=\varepsilon\\
h&\mapsto h^*+(gh)^*\\
g&\mapsto 1^*-g^*\\
gh&\mapsto  h^*-(gh)^*.
\end{array}$$
So, the cotriangular structures of $H_4$ can be obtained applying the isomorphism $\phi\otimes\phi$ to the $R_t$'s. They are:
$$\begin{array}{c|rrrr}
r_t&1&g&h&gh \\
\hline
1&1&1&0&0\\
g&1&-1&0&0\\
h&0&0&t&-t\\
gh&0&0&t&t\\
\end{array}$$

The Drinfeld double $D(H_4)=H_4^{*,cop}\bowtie H_4$ of $H_4$ is
isomorphic to the Hopf algebra generated by $\phi(h)\bowtie 1$, $\phi(g)\bowtie
1$, $\varepsilon\bowtie g$ and $\varepsilon\bowtie h$
with relations:
$$\begin{array}{l}
(\phi(h)\bowtie 1)^2=0;\\
(\phi(g)\bowtie 1)^2=\varepsilon\bowtie 1;\\
(\phi(h)\bowtie 1)(\phi(g)\bowtie 1)+(\phi(g)\bowtie 1)(\phi(h)\bowtie 1)=0;\\
(\varepsilon\bowtie h)^2=0; \\
(\varepsilon\bowtie h)(\varepsilon\bowtie g)+(\varepsilon\bowtie g)(\varepsilon\bowtie h)=0;\\
(\varepsilon\bowtie g)^2=\varepsilon\bowtie 1;\\
(\phi(h)\bowtie 1)(\varepsilon\bowtie g)+(\varepsilon\bowtie
  g)(\phi(h)\bowtie 1)=0;\\
(\phi(g)\bowtie 1)(\varepsilon\bowtie  h)+(\varepsilon\bowtie h)(\phi(g)\bowtie 1)=0;\\
(\varepsilon\bowtie  g)(\phi(g)\bowtie 1)=(\phi(g)\bowtie 1)(\varepsilon\bowtie g);\\
(\phi(h)\bowtie 1)(\varepsilon\bowtie h)-(\varepsilon\bowtie h)(\phi(h)\bowtie 1)=(\phi(g)\bowtie 1)-(\varepsilon\bowtie g)
\end{array}$$
and with coproduct induced by the coproducts in $H_4$ and $H_4^{*,cop}$. For $l\in H_4$ we will sometimes write $\phi(l)$ instead of $\phi(l)\bowtie 1$ and $l$ instead of $1\bowtie l$ for simplicity. \smallskip

Let us recall that a Yetter-Drinfeld $H_4$-module $M$ with action $\cdot$ and coaction $\rho$
becomes a $D(H_4)$-module by letting $1\bowtie l$ act as $l$ for every
$l\in H_4$ and $(\phi(l)\bowtie 1).m=\m0(\phi(l)(\m1))$ for $m \in M$. Conversely, a
$D(H_4)$-module $M$ becomes naturally a Yetter-Drinfeld module with
$H_4$-action obtained by restriction and $H_4$-coaction given by
$$\rho(m)=\frac{1}{2}(\phi(1+g).m\otimes 1+\phi(1-g).m\otimes
g+\phi(h+gh).m\otimes h+\phi(h-gh)\otimes gh).$$

We will often switch from one notation to the other according to convenience. \par \bigskip

{\it Centers and centralizers.} If $A$ is a Yetter-Drinfeld $H$-module algebra,
and $B$ is a Yetter-Drinfeld submodule algebra of $A$, the left and
the right centralizer of $B$ in $A$ are defined to be:
$$C^l_A(B):=\{a\in A~|~ba=\a0(\a1\cdot b)\ \forall b\in B\},$$
$$C^r_A(B):=\{a\in A~|~ab=\b0(\b1\cdot a)\ \forall b\in B\}.$$
For the particular case $B=A$ we have the right center $Z^r(A)$ and the left center
$Z^l(A)$ of $A$. Both are trivial when $A$ is $H$-Azumaya, \cite[Proposition 2.12]{CVZ2}.

\section{Some low dimensional representatives in $BQ(k, H_4)$}

In this section we shall introduce a family of 2-dimensional representatives of classes in $BQ(k,H_4)$ that will turn out to be
easy to compute with. They appeared for the first time in \cite{MS} and a particular case of them is treated in
\cite[Section 1.5]{ACZ}. \smallskip

Let $a,\,t,\,s\in k$. The algebra $C(a)$ generated by $x$ with relation $x^2=a$ is acted upon by $H_4$ by
$$g\cdot 1=1,\quad g\cdot x=-x,\quad h\cdot 1=0, \qquad h\cdot x=t,$$
and it is a right $H_4$-comodule via
$$\rho_s(1)=1\otimes 1,\quad\quad\rho_s(x)=x\otimes g+s\otimes h.$$

It is not hard to check that $C(a)$ with this action and coaction  is a left $H_4$-module algebra and a right
$H^{op}$-comodule algebra. We shall denote it by $C(a;t,s)$.

\begin{lemma}\label{properties} Let notation be as above.
\begin{enumerate}
\item[(1)] $C(a;t,s)$ is a Yetter-Drinfeld module algebra with the preceding structures.

\item[(2)] As a module algebra  $C(a;t,s)\cong C(a';t',s')$ if and only if
there is $\alpha\in k^\cdot$ such that $a=\alpha^2a'$ and $t=\alpha t'$.

\item[(3)] As a comodule algebra  $C(a;t,s)\cong C(a';t',s')$ if and only if
 there is $\alpha\in k^\cdot$ such that $a=\alpha^2a'$ and $s=\alpha s'$.

\item[(4)] As a Yetter-Drinfeld module algebra $C(a;t,s)\cong C(a';t',s')$ if and only if
there exists $\alpha\in k^\cdot$ such that $a=\alpha^2a'$, $t=\alpha t'$ and $s=\alpha s'$.

\item[(5)] The module structure on $C(a;t,s)$ is induced by its comodule structure and a
cotriangular structure $r_l$ if and only if $t=sl$.

\item[(6)] The comodule structure on $C(a;t,s)$ is induced by its module
  structure and a triangular structure $R_l$ if and only if $s=lt$.

\item[(7)] The $H_4$-opposite algebra of $C(a;t,s)$ is $C(st-a;t,s)$.

\item[(8)] $C(a;t,s)$ is an $H_4$-Azumaya algebra if and only if $2a \neq st$.
\end{enumerate}
\end{lemma}

\noindent \pf Let $x$ and $y$ be algebra generators in $C(a;t,s)$ and $C(a';t',s')$ respectively with $x^2=a$ and $y^2=a'$. \smallskip

\noindent (1) We verify condition (\ref{YD}) for $b=x$ and $l=h$. The other cases are easier to check.
$$\begin{array}{l}
\h2\cdot x_{(0)} \otimes \h3 x_{(1)} S^{-1}(\h1) \\
\hspace{1.7cm} = g\cdot x\otimes(-gh)+g\cdot s\otimes (gh)(-gh) +h\cdot x\otimes   g^2 \\
\hspace{2cm} +h\cdot s\otimes gh+x\otimes hg+s\otimes h^2 \\
\hspace{1.7cm} = x\otimes gh+t\otimes 1-x\otimes gh \\
\hspace{1.7cm} = \rho_s(h\cdot x).
\end{array}
$$

\noindent (2) An algebra isomorphism $f\colon C(a;t,s)\to C(a';t',s')$ must map $x$ to $\alpha
y$ for some $\alpha\in k^\cdot$. Then $a=x^2=(\alpha y)^2=\alpha^2 a'$. Besides, $h.f(x)=f(h.x)$ implies $t'\alpha=t$. It is easy
to verify that the condition is also sufficient. \vspace{5pt}

\noindent (3) In the above setup $\rho_{s'}(f(x))=(f\otimes\id)\rho_s(x)$ implies $s'\alpha=s$. It is
not hard to check that this condition is also sufficient. \vspace{5pt}

\noindent (4) It follows from the preceding statements. \vspace{5pt}

\noindent (5) If the module structure on $C(a;t,s)$ is induced by its comodule structure $\rho_s$ and
some $r_l\in{\cal U},$ then $t=h\cdot x=xr_l(h\otimes g)+sr_l(h\otimes h)=sl.$ Conversely, if
$t=sl$, then
$$\begin{array}{l}
g\cdot 1=1=1r_l(g\otimes 1); \qquad h\cdot 1=0=1r_l(h\otimes 1); \vspace{2pt} \\
g\cdot x=-x=xr_l(g\otimes g)+sr_l(g\otimes h)=x_{(0)}r_l(g\otimes x_{(1)}); \vspace{2pt} \\
h\cdot x=t = xr_l(h \otimes g)+sr_l(h \otimes h)=x_{(0)}r_l(h\otimes x_{(1)}).
\end{array}$$
Therefore the action is induced by the coaction and $r_l$. \vspace{5pt}

\noindent (6)  If the comodule structure on $C(a;t,s)$ is induced by the
  action and some $R_l\in{\cal T},$ then
$$x\otimes g+s\otimes h=\rho_s(x)=(R_l^{(2)}\cdot x)\otimes R_l^{(1)}=\frac{1}{2}(2x\otimes g)+\frac{l}{2}(2t\otimes h)=x\otimes g+lt\otimes h
$$
hence $s=lt$. Conversely, if $s=lt$ then
$$\begin{array}{l}
\rho_s(1)=1\otimes 1=(R_l^{(2)}\cdot 1)\otimes R_l^{(1)}, \\
\rho_s(x)=x\otimes g +s \otimes h= (R_l^{(2)}\cdot x)\otimes R_l^{(1)},
\end{array}$$
so the comodule structure is induced by the action and $R_l$. \vspace{5pt}

\noindent (7) $\overline{C(a;t,s)}$ has $1,\,x$ as a basis and $1$ is the unit.
The action and coaction on $1$ and $x$ are as for $C(a;t,s)$. By direct computation,
$x\circ x=x(g\cdot x)+s(h\cdot x)=-a+st,$ so
$\overline{C(a;t,s)}=C(st-a;t,s)$. \vspace{5pt}

\noindent (8) The algebra $C(a;t,s)$ is $H_4$-Azumaya if and only if the maps $F$ and $G$
defined in \eqref{f-and-g} are isomorphisms. The space
$C(a;t,s)\# C(a;t,s)$ has ordered basis $1\#1,\,1\# x,\,x\#1,\,x\# x$ while ${\rm End}(C(a;t,s))$ has
basis $1^*\otimes 1,\,1^*\otimes x,\,x^*\otimes 1,\,x^*\otimes x$ with
the usual identification $C(a;t,s)^*\otimes C(a;t,s)\cong{\rm End}(C(a;t,s))$.
Then for every $b,c \in C(a;t,s)$ we have
$$\begin{array}{l}
F(b\#c)(1)=bc,\quad F(b\# c)(x)=bx(g\cdot c)+sb(h\cdot c), \vspace{2pt} \\
G(1\# b)(c)=cb,\quad G(x\# b)(c)=x(g\cdot c)b+s(h\cdot c)b.
\end{array}$$
The matrices associated with $F$ and $G$ with respect to the given bases are respectively
$$\left(
\begin{array}{cccc}
1&0&0&a\\
0&1&1&0\\
0&st-a&a&0\\
1&0&0&st-a
\end{array}
\right) \qquad \left(
\begin{array}{cccc}
1&0&0&a\\
0&1&1&0\\
0&a&st-a&0\\
1&0&0&st-a
\end{array}
\right)$$
whose determinants $-(st-2a)^2$ and $(st-2a)^2$ are nonzero if and only if $2a\neq st$. \hfill $\Box$ \medskip

We have seen so far that the algebras $C(a;s,t)$ can be viewed as representatives of classes in $BM(k,H_4, R_l)$ or in $BC(k,H_4,r_l)$ for suitable $l \in k$. It is known that these groups are all isomorphic to $(k,+) \times BW(k),$ where $BW(k)$ is the Brauer-Wall group of $k$. We aim to find to which pair $(\beta, [A])\in (k,+) \times BW(k)$ do the class of $C(a;t,s)$ correspond. The group $BM(k,H_4,R_0)$  was computed in \cite{VOZ3}. The computation of $BC(k,H_4,r_0)$ follows from self-duality of $H_4$. It was shown in \cite{gio} that all groups $BC(k,H_4,r_t)$ (hence, dually, all $BM(k,H_4,R_t)$) are isomorphic. We shall use the description of $BM(k,E(1),R_t)$ given in \cite{enne} beause this might allow generalizations. In the mentioned paper the Brauer group $BM(k,E(n),R_0)$ is computed for the family of Hopf algebras $E(n)$, where $E(1)=H_4$. We shall recall first where do the isomorphism of the different Brauer groups $BC$ and $BM$ stem from. The notion of lazy cocycle plays a key role here.

We recall from \cite{lazy} that a lazy cocycle on $H$ is a left 2-cocycle $\sigma$ such that twisting $H$ by $\sigma$ does
not modify the product in $H$. In other words: for every $h, l, m\in H,$
\begin{equation}
\sigma(\h1\otimes\l1)\sigma(\h2\l2\otimes
m)=\sigma(\l1\otimes\m1)\sigma(h\otimes \l2\m2)
\end{equation}
\begin{equation}
\sigma(\h1\otimes\l1)\h2\l2=\h1\l1\sigma(\h2\otimes\l2)
\end{equation}

It turns out that a lazy left cocycle is also a right cocycle. Given a
lazy cocycle $\sigma$ for $H$ and a $H^{op}$-comodule algebra $A$, we may
construct a new $H^{op}$-comodule algebra $A_{\sigma}$, which is equal to $A$ as
a $H^{op}$-comodule, but with product defined by:
$$a\bullet b=\a0\b0\sigma(\a1\otimes\b1).$$

The group of lazy cocycles for $H_4$ is computed in
\cite{lazy}. Lazy cocycles are parametrized by elements $t \in k$
as follows:
$$\begin{array}{c|rrrr}
\sigma_t&1&g&h&gh\\
\hline
1&1&1&0&0\\
g&1&1&0&0\\
h&0&0&\frac{t}{2}&\frac{t}{2}\\
gh&0&0&\frac{t}{2}&-\frac{t}{2}\\
\end{array}
$$

We have the following group isomorphisms:
\begin{enumerate}
\item[(2.3)] $\Psi_t: BC(k,H_4,r_0) \rightarrow BC(k,H_4,r_t), [A] \mapsto[A_{\sigma_t}],$ constructed in \cite[Proposition 3.1]{gio}.

\item[(2.4)] $\Phi_t:BM(k,H_4,R_t) \rightarrow BC(k,H_4,r_t), [A] \mapsto [A^{op}].$ We explain how $A^{op}$ is equipped with the corresponding structure. The left $H_4$-module algebra $A$ becomes a right $H_4^*$-comodule algebra. Then $A^{op}$ is a right $H_4^{*,op}$-comodule algebra. The quasitriangular structure $R_t$ is a coquasitriangular structure in $H_4^*$. Then $A^{op}$ may be endowed with the left $H_4^*$-action stemming from the comodule structure and $R_t$. On the other hand, $A^{op}$ may be viewed as an $H_4^{op}$-comodule algebra through the isomorphism $\phi:H_4 \rightarrow H_4^*$. The coquasitriangular structure $R_t$ on $H_4^*$ corresponds to the coquasitriangular structure $r_t$ on $H_4$ via $\phi.$
\end{enumerate}
\setcounter{equation}{4}

An isomorphism between $BM(k,H_4,R_0)$ and $BM(k,H_4,R_t)$ can be constructed combining the above ones. Thus, the crucial step is to analyze the sought correspondence for $BM(k,H_4,R_0)$. \smallskip

The Brauer group $BM(k, H_4, R_0)$ is computed in \cite{VOZ3} through the split exact
sequence (see also \cite[Theorem 3.8]{ACZ} for an alternative approach):
$$\xymatrix{1\ar[r] & (k,+) \ar[r] & BM(k, H_4,R_0) \ar@<0.5ex>[r]^(0.6){j^*} &
BW(k) \ar@<0.5ex>[l]^(0.4){\pi^*} \ar[r] & 1.}$$
The map $j^*:BM(k, H_4,R_0)\to BW(k),[A] \mapsto [A]$ is obtained by restricting the
$H_4$-action of $A$ to a $k{\mathbb Z}_2$-action via the inclusion map
$j:k{\mathbb Z}_2 \rightarrow H_4$. This map is split by $\pi^*:BW(k) \rightarrow BM(k, H_4,R_0),
[B] \mapsto [B]$, where $B$ is considered as an $H_4$-module by restriction of scalars via
the algebra projection $\pi:H_4 \rightarrow k{\mathbb Z}_2, g\mapsto g, h \mapsto 0$. A class $[A]$ lying in the
kernel of $j^*$ is a matrix algebra with an inner action of $H_4$ such that the restriction to $k{\mathbb Z}_2$ is
strongly inner. Thus there exist uniquely determined $u, w \in A$ such that
\begin{equation}\label{inv1}
g\cdot a=uau^{-1}, \quad h\cdot a= w(g\cdot a)-aw \quad\forall a\in A,
\end{equation}
\begin{equation}
\label{inv2}u^2=1, \quad wu+uw=0, \quad w^2=\beta,
\end{equation}
for certain $\beta \in k$. Mapping $[A]\mapsto \beta$ defines a group isomorphism $\chi\colon Ker(j^*) \\ \cong (k, +)$. We will determine $j^*([C(a;t,s)])$ and $\chi([C(a;t,s)]\pi^*j^*([C(a;t,s)]^{-1}))$ whenever this is well-defined. To this purpose, we will first describe all products of two algebras of type $C(a;t,s)$.

\begin{lemma}\label{product} Let $x,y$ be generators for $C(a;t,s)$ and $C(a';t',s')$ respectively, with relations, $H_4$-actions and coactions as above. The product $C(a;t,s)\# C(a';t',s')$ is isomorphic to the generalized quaternion algebra
with generators $X=x\#1$ and $Y=1\# y$, relations, $H_4$-action and and $H_4$-coaction:
$$X^2=a,\quad Y^2=a', \quad XY+YX=st',$$
$$g\cdot X=-X,\quad g\cdot Y=-Y,\quad h\cdot X=t,\quad h\cdot Y=t',$$
$$\rho(X)=X\otimes g+s\otimes h,\quad \rho(Y)=Y\otimes g+s'\otimes h.$$
\end{lemma}

\noindent \pf By direct computation:
$$X^2=(x\#1)(x\#1)=a\# 1,\quad Y^2=(1\#y)(1\#y)=a'\#1,\quad XY=x\# y,$$
$$YX=(1\# y)(x\#1)=x\#(g\cdot y)+s\#(h\cdot y)=-XY+st'\#1.$$
The formulas for the action and the coaction follow immediately from the
definition of action and coaction on a $\#$-product. \hfill $\Box$ \smallskip

Elements in $BW(k)$ are represented by graded tensor products of the following three type of algebras:
$C(1)$ generated by the odd element $x$ with $x^2=1$; classically Azumaya algebras having trivial ${\mathbb
  Z}_2$-action; and $C(a)\# C(1),$ where $C(a)$ is generated by the
odd element $y$ with $y^2=a\in k^\cdot$ (\cite[Theorem IV.4.4]{lam}).

\begin{proposition}\label{BM0} For $a\neq0$ let $[C(a;t,0)] \in BM(k, H_4,R_0)$ denote the class of $C(a;t,0).$ Then $$[C(a;t,0)]=(t^2(4a)^{-1}, [C(a)])\in (k,+)\times BW(k),$$
so the group $BM(k, H_4, R_0)$ is generated by $BW(k)$ and the
classes $[C(a;1,0)]$.
\end{proposition}

\noindent \pf It is clear that if $a\neq0$ then $j^*([C(a;t,0)])=[C(a)]$
and that $\pi^*([C(a)])=[C(a;0,0)]$. Thus, $[C(a;t,0)\#
  C(-a;0,0)]\in{\rm Ker}(j^*)$. We shall compute its image through $\chi$. By Lemma \ref{product}, $C(a;t,0)\# C(-a;0,0)$
is generated by $X$ and $Y$ with relations, $H_4$-action and $H_4$-coaction:
$$X^2=a,\quad Y^2=-a,\quad XY+YX=0,$$
$$g\cdot X=-X,\quad g\cdot Y=-Y,\quad h\cdot X=t,\quad h\cdot Y=0,$$
$$\rho(X)=X\otimes g,\quad \rho(Y)=Y\otimes g.$$
We look for the element $w$ satisfying (\ref{inv1}) and
(\ref{inv2}). This element must be odd with respect to the ${\mathbb
  Z}_2$-grading induced by the $g$-action, hence $w=\lambda X+\mu Y$
for some $\lambda,\mu\in k$. Condition $h\cdot X=-wX-Xw$ implies $t=-2\lambda a$ and
condition $h\cdot Y=-wY-Yw$ implies $0=-2\mu a$ so $w^2=a\lambda^2=t^2(4a)^{-1}$.
Thus $[C(a;t,0)]=(t^2(4a)^{-1}, [C(a)])$ and we have the first statement. For the second one,
let $(\beta, [A])\in (k,+)\times BW(k)$. If $\beta=0$ there is nothing to prove. If $\beta\neq 0$, the class
$[C((4\beta)^{-1}t^2;t,0)]=[C((4\beta)^{-1};1,0)]=(\beta,
[C((4\beta)^{-1})])$, so $BM(k,H_4,R_0)\cong(k,+)\times BW(k)$ is
generated by $BW(k)$ and the $[C(a;1,0)]$ for $a\neq0.$\hfill$\Box$

\begin{lemma}\label{BWinBM} Let $A$ be a $D(H_4)$-module algebra.
\begin{enumerate}
\item[(1)] If the $h$-action on $A$ is trivial, then $A$ is $(H_4,R_0)$-Azumaya if and only if it is $(H_4,R_t)$-Azumaya for every $t\in k$.
\item[(2)] If the $\phi(h)$-action on $A$ is trivial, then $A$ is $(H_4,r_0)$-Azumaya if and only if it is $(H_4,r_t)$-Azumaya for every $t\in k$.
\item[(3)] The representatives of $BW(k)$ in $BC(k, H_4, r_t)$ and $BM(k, H_4, R_s)$ all coincide when viewed
  inside $BQ(k,H_4)$.
\end{enumerate}
\end{lemma}

\noindent \pf (1) It follows from the form of the elements in ${\cal T}$ that if $A$ is $(H_4,R_0)$-Azumaya
and the action of $h$ on $A$ is trivial (i.e., if it lies in $BW(k)$), then
its comodule structure $\rho_t$ induced by $R_t$ coincides with
the comodule structure $\rho_0$ induced by $R_0$. Hence, the maps $F$
and $G$ with respect to the action and $\rho_t$ are the same as the
maps $F$ and $G$ with respect to the action and $\rho_0$, so $A$ is
$(H_4,R_t)$-Azumaya for every $t\in k$. \vspace{5pt}

(2) It is proved  as (1). \vspace{5pt}

(3) The first statement shows that the representatives of $BW(k)$ inside the different $BM(k, H_4, R_t)$
coincide. The second statement shows the same for $BC(k, H_4, r_t)$. Therefore we may assume $s=t=0$.
 The elements of this copy of $BW(k)$ consist of ${\mathbb Z}_2$-graded Azumaya algebras $A$ where the grading
 is induced by the action of $g$. The $h$-action is trivial. If the coaction $\rho$ is induced by $R_0$, then $a \in A$ is
 odd if and only if $\rho(a)=a\otimes g$. The action $\ru$ induced on $A$ by $r_0$ and $\rho$ is as follows:
$h\ru a=0$ for every $a \in A$ and $g\ru a=-a$ if and only if $\rho(a)=a\otimes g$, that is, the original action on $A$ and $\ru$
coincide. Thus, the maps $F$ and $G$ coincide in all cases and $A$ represents an element in $BW(k)\subset BM(k,H_4,R_0)$ if and only if it represents an element in $BW(k)\subset BC(k,H_4,r_0)$. \hfill$\Box$

\begin{proposition}\label{BCs}
The group $BC(k, H_4, r_s)$ is generated by the Brauer-Wall group and the classes $[C(a;s,1)]$ for $2a \neq s$.
\end{proposition}

\noindent \pf We will first deal with the case $s=0$. We will show that the isomorphism $\Phi_0: BM(k, H_4,R_0) \rightarrow BC(k,H_4,r_0),$ $[A] \mapsto [A^{op}]$ in (2.4) maps  $[C(a;1,0)]$ to $[C(a;0,1)]$ and $BW(k)\subset BM(k, H_4,R_0)$ to $BW(k)\subset BC(k, H_4,r_0)$. The class $[C(a;1,0)]$ is mapped to the class of the algebra $C(a)^{op}$ with comodule
structure
$$\rho(x)=x\otimes(1^*-g^*)+1\otimes (h^*+(gh)^*)=x\otimes\phi(g)+1\otimes\phi(h)$$
and $H_4$-action induced by the cotriangular structure $r_0$, that is, $g \cdot x=-x$ and $h\cdot x=0$. The algebra $C(a)^{op}$ with these structures is just $C(a;0,1)$.

Let $A$ be a representative of a class in $BW(k)\subset BM(k, H_4,R_0)$ with action $\cdot $ for which $h\cdot a=0$ for all $a \in A$. The class $ [A]$ is mapped by $\Phi_0$ to the class of $A^{op}$ with coaction
$$\rho(a)=a\otimes 1^*+(g\cdot a)\otimes g^*+(h\cdot a)\otimes h^*+(gh\cdot a)\otimes (gh)^*\in A\otimes \phi(k{\mathbb Z}_2).$$
Therefore $[A^{op}] \in BW(k) \subset BC(k, H_4,r_0)$.

We now take $s \in k$ arbitrary and use the isomorphism $\Psi_s: BC(k,H_4,r_0) \rightarrow BC(k,H_4,r_s)$ in (2.3) to prove the statement. We will show that $[C(a;0,1)]$ is mapped to $[C(a+2^{-1}s;s,1)]$ through $\Psi_s$. Recall that $\Psi_s$ maps the class of $C(a;0,1)$ to the class of the algebra
$C(a;0,1)_{\sigma_s}$. It is generated by $x$ with relation
$$x\bullet x=x^2\sigma_s(g\otimes g)+x\sigma_s(h\otimes
g)+x\sigma_s(g\otimes h)+\sigma_s(h\otimes h)=a+\frac{s}{2},$$
with (same) coaction $\rho(x)=x\otimes g+1\otimes h$ and action induced by $\rho$ and
$r_s$, that is:
$$g\cdot x=r_s(g\otimes g)x+r_s(g\otimes h)=-x, \quad h\cdot x=r_s(h\otimes g)x+r_s(h\otimes h)=s.$$
Then $\Psi_s([C(a;0,1)])=[C(a+\frac{s}{2}; s,1)]$.

Since the coaction is not changed by $\Psi_s$ the class of an element $A$ for which the image of the coaction is in $A \otimes k{\mathbb Z}_2$ is again of this form. Hence the classes in $BW(k)\subset BC(k, H_4, r_0)$ correspond to the classes in $BW(k) \subset BC(k, H_4, r_s)$. \hfill $\Box$

\begin{proposition}\label{BMt}
The group $BM(k, H_4, R_t)$ is generated by the Brauer-Wall group and the classes $[C(a;1,t)]$ for $2a \neq t$.
\end{proposition}

\noindent \pf Through the isomorphism  $\Phi_t:BM(k, H_4, R_t) \rightarrow BC(k,H_4,r_t)$ in (2.4), the class $[C(a;1,t)]$ is mapped to $[C(a;t,1)]$ and the classes in $BW(k)\subset BM(k, H_4, R_t)$ correspond to the classes in $BW(k)\subset BC(k, H_4, r_t)$. The $H_4$-comodule structure on the algebra $C(a)^{op}$ is:
$$\rho(x)=x\otimes(1^*-g^*)+1\otimes (h^*+(gh)^*)=x\otimes\phi(g)+1\otimes\phi(h)$$
The $H_4$-action induced by the cotriangular structure $r_t$ on $H_4$ gives $h\cdot x=t$. Therefore this algebra is $C(a;t,1)$. Finally, the statement concerning $BW(k)$ is proved as in the preceding theorem.\hfill$\Box$

\begin{remark}{\em That $BM(k,H_4, R_t)$ is generated by $BW(k)$ and the classes $[C(a;1,t)]$ for $2a \neq t$ was first discovered in \cite[Theorem 3.8 and Page 392]{ACZ} as a consequence of the Structure Theorems for $(H_4,R_t)$-Azumaya algebras. Since we will strongly use Proposition \ref{BMt} later, for the reader's convenience we offered this alternative and self-contained approach. Notice that it mainly relies on Lemma \ref{product} that will be another key result for us in the sequel.}
\end{remark}

\section{Fitting $BM(k, H_4,R_t)$ and $BC(k, H_4, r_s)$ into $BQ(k, H_4)$}

As groups $BM(k,H_4,R_t)\cong BC(k, H_4, r_s)$ for every $s,t \in k$. However, their images in $BQ(k,H_4)$ through
the natural embeddings
$$i_t\colon BM(k, H_4, R_t)\to BQ(k, H_4) \quad \textrm{and} \quad
\iota_s\colon BC(k, H_4, r_s)\to BQ(k, H_4)$$ do not coincide in general. In
this section we will describe the mutual intersections of these images.

\begin{proposition}\label{BM=BC}
Let $0\neq t \in k$ then $Im(i_t)=Im(\iota_{t^{-1}})$
\end{proposition}

\noindent \pf Given $t \neq 0$, by Lemma \ref{properties},  $[C(a;1,t)] \in Im(i_t) \cap Im(\iota_{t^{-1}})$
for every $a\neq 2t.$ Besides, by Lemma \ref{BWinBM},
$i_t(BW(k))=\iota_s(BW(k))$ for any $s \in k$. Since the elements of $BW(k)$ and
the $[C(a;1,t)]$'s generate $BM(k, H_4, R_t)$ and $BC(k, H_4, r_{t^{-1}})$
we are done.\hfill$\Box$ \medskip

Given $[A]$ in $BQ(k,H_4)$, there are two natural ${\mathbb Z}_2$-gradings on $A$, the one coming from
the $g$-action, for which $|a|=1$ iff $g\cdot a=-a$ for $0 \neq a \in A$ and the one arising from
the coaction, for which $\deg(a)=1$ if and only if
$(\id\otimes\pi)\rho(a)=a\otimes g$ where $\pi$ is the projection onto
$k{\mathbb Z}_2$. If we view $A$ as a
$D(H_4)$-module, the grading $|\cdot |$ is associated with the $1\bowtie g$-action whereas the grading
$\deg$  is associated with the $\phi(g)\bowtie 1$-action. Let us observe that for the classes
$C(a;t,s)$ the two natural gradings coincide, for every $a,\,t,\,s \in k$.

\begin{lemma}\label{grading}
Let $[A]\in BQ(k,H_4)$ and $[B]$ in $i_0(BW(k))$. As a $H_4$-module algebra,
\begin{itemize}
\item[(1)] $A\# B\cong A \widehat{\otimes} B$, the ${\mathbb Z}_2$-graded tensor product with respect to the
$\deg$-grading on $A$ and the natural $|\cdot|$-grading on $B$.

\item[(2)] $B\# A\cong B \widehat{\otimes} A$, the ${\mathbb Z}_2$-graded tensor product with respect to the $|\cdot |$-grading on $A$ and the natural $|\cdot |$-grading on $B$.
\end{itemize}
\end{lemma}

\noindent \pf The two gradings on $B$ coincide and we have, for homogeneous $b\in B$
and $c\in A$ (for the $\deg$-grading):
$$(a\# b)(c\# d)=a\c0\# (\c1\cdot b)d=a c\# (g^{\deg(c)}\cdot b)d=(-1)^{\deg(c)|b|}ac\# bd.$$
For homogeneous $b\in B$ and $c \in A$ (for the $|\cdot |$-grading):
$$(d\# c)(b\# a)=d\b0\# (\b1\cdot c)a=db\# (g^{|b|}\cdot c)a=(-1)^{|c||b|}db\# ca.$$
\hfill$\Box$

It follows from Propositions \ref{BCs}, \ref{BMt} and Lemma \ref{grading} that all elements in $Im(i_t)$ and $Im(\iota_t)$ can be represented by algebras for which the two ${\mathbb Z}_2$-gradings coincide, since this property is respected by
the $\#$-product. Indeed, this kind of representatives give rise to a subgroup that we will study in Section 5. \smallskip

We will show now that groups of type $BC$ or $BM$ either intersect
only in $BW(k)$ or coincide and that the latter happens only in the
situation of Proposition \ref{BM=BC}.

\begin{theorem}\label{distinct}
Consider the class of $C(a;t,s)$ in $BQ(k,H_4)$. Then:
\begin{enumerate}
\item[(1)] $[C(a;t,s)]\in Im(i_l)$ if and only if $s=lt$;
\item[(2)]  $[C(a;t,s)]\in Im(\iota_l)$ if and only if $sl=t$.
\end{enumerate}
\end{theorem}

\noindent \pf (1) We know from Lemma \ref{properties} that if the action (resp. coaction) of $C(a;t,s)$ comes from the cotriangular (resp. triangular) structure, then the indicated relations among the parameters hold. We only
need to show that the condition is still necessary if we change representative in the class.

Let us assume that $[C(a;t,s)]\in Im(i_l)$ for some $l \in k$. Then $[C(a;t,s)]=[C(b;1,l)][A]=[A][C(b;1,l)]$ for some $[A]\in
i_l(BW(k))$ and $b \in k$ with $2b \neq l$. Hence $[C(a;t,s)\# C(l-b;1,l)]=[A]\in i_l(BW(k))$. We
may choose $A$ so that the $h$-action and the $\phi(h)$-action on $A$ are trivial.

Since $[C(a;t,s)\# C(l-b;1,l)\#\overline{A}]$ is trivial in $BQ(k,H_4)$, there is a $D(H_4)$-module $P$ such that
$C(a;t,s)\# C(l-b;1,l)\# \overline{A}\cong{\rm End(P)}$ as $D(H_4)$-module algebras. Then
${\rm End(P)}$ has a strongly inner $D(H_4)$-action. In other words, there is a convolution invertible
algebra map $\nu\colon D(H_4)\to {\rm End}(P)$ such that
$$(m\bowtie n)\cdot f=\nu(\m2\bowtie\n1) f \nu^{-1}(\m1\bowtie\n2)$$
for every $m\bowtie n\in D(H_4), f\in {\rm End(P)}$, where $\nu^{-1}$ denotes the convolution
inverse of $\nu$. In particular, for $u=\nu(\varepsilon\bowtie g)$ and $w=\nu(\varepsilon\bowtie h)u$ we have
$$\begin{array}{l}
g\cdot f=u f u^{-1}, \quad h\cdot f=w(g\cdot f)-f w, \vspace{5pt} \\
u^2=1, \quad w^2=0, \quad u w+w u=0.
\end{array}$$

We should be able to find $U,\,W\in C(a;t,s)\# C(l-b;1,l)\#\overline{A}$ such that
$$\begin{array}{l}
U^2=1, \quad g\cdot Z=UZU^{-1}, \vspace{5pt} \\
g\cdot W=-W, \quad W^2=0, \quad h\cdot Z=W(g\cdot Z)-ZW
\end{array}$$
for all $Z$ in $C(a;t,s)\# C(l-b;1,l)\#\overline{A}.$

Using the presentation of $C(a;t,s)\# C(l-b;1,l)$ in Lemma \ref{product} we may
write $W=\sum_{0\le i,j\le 1}X^iY^j\#\alpha_{ij}$ with $\alpha_{ij}\in\overline{A}$
homogeneous of degree $i+j+1\;{\rm mod}\ 2$ with respect to the $g$-grading.
Since the action of $h$ on $1\# \overline{A}$ is trivial we have, for
homogeneous $\gamma\in\overline{A}$:
$$\begin{array}{ll}
0 &=h\cdot (1\# \gamma) \vspace{2pt}\\
   &=W(g\cdot (1\#\gamma))-(1\# \gamma)W  \vspace{2pt} \\
   &=(-1)^{|\gamma|}\sum_{0\le i,j\le 1}X^iY^j\#\alpha_{ij}\gamma-\sum_{0\le i,j\le 1}(X^iY^j)_{(0)}\#((X^iY^j)_{(1)}\cdot  \gamma)\alpha_{ij} \vspace{2pt} \\
&=(-1)^{|\gamma|}[1\#\alpha_{00}\gamma+Y\#\alpha_{01}\gamma+X\#\alpha_{10}\gamma+XY\#\alpha_{11}\gamma] \vspace{2pt} \\
&\phantom{=}-1\#\gamma\alpha_{00}-Y\#(-1)^{|\gamma|}\gamma\alpha_{01}-X\#(-1)^{|\gamma|}\gamma\alpha_{10}-XY\#\gamma\alpha_{11}.
\end{array}$$
From here we deduce that the odd elements $\alpha_{00}, \alpha_{11}$ and the even elements $\alpha_{10}, \alpha_{01}$ belong to the ${\mathbb Z}_2$-center of $\overline{A}.$ Hence $\alpha_{00}, \alpha_{11}$ are zero and $\alpha_{10}, \alpha_{01}$ are scalars. So, we can write $W=\alpha X\# 1+\beta Y\#1$  for some $\alpha,\beta\in k$ and we will get:
$$\begin{array}{rl}
\alpha t+\beta & =h\cdot W=-2W^2=0,\\
t & =h \cdot (X\# 1)=\alpha(-2a+t s),\\
1 & =h \cdot (Y\# 1)=-\alpha s-2\beta (l-b)=\alpha(-s+2t(l-b)).
\end{array}$$
Combining the second equation with the third one multiplied by $t$ and using $\alpha\neq 0$ we obtain
\begin{equation}\label{a=}
a=ts-t^2(l-b).
\end{equation}

The $|\cdot |$-grading and the $\deg$-grading on $C(a;t,s)\#C(l-b;1,l)\# \overline{A}$ coincide. Therefore:
$$\nu(\phi(g)\bowtie 1)f\nu(\phi(g)\bowtie 1)^{-1}=\phi(g)\cdot f= g\cdot f=u f u^{-1} \qquad \forall f \in {\rm End}(P).$$
Since ${\rm End}(P)$ is central and $\nu$ is an algebra morphism,
$u':=\nu(\phi(g)\bowtie 1) =\lambda u$ with $\lambda=\pm 1$ (both possibilities will be analyzed later). The element $w':=\nu(\phi(h)\bowtie 1)$ satisfies
$$\phi(h)\cdot f =w'f-(\phi(g) \cdot f)w' \qquad \forall f \in {\rm End}(P).$$
Thus, we can take $W'$ in $C(a;t,s)\# C(l-b;1,l)\# \overline{A}$ such that
$$W'U+UW'=0, \quad (W')^2=0 \quad \phi(h)\cdot Z=W'Z-(g\cdot Z)W'$$
for all $Z$ in $C(a;t,s)\# C(l-b;1,l)\#\overline{A}.$ Arguing as for $W$ before, we see that $W'=\gamma X\# 1+\delta Y\# 1$ for
some $\gamma, \delta\in k$. It follows from the last relation of $D(H_4)$ in \S \ref{preliminaries} that
$$\nu(\varepsilon\bowtie hg)\nu(\phi(h)\bowtie 1)+\nu(\phi(h)\bowtie 1)\nu(\varepsilon\bowtie hg)=\nu(\phi(g)\bowtie 1)\nu(\varepsilon\bowtie g)-\nu(\varepsilon\bowtie g)^2.$$
This implies $WW'+W'W=\lambda-1.$ Besides,
$$0=\phi(h)\cdot W'=2(W')^2=s\gamma+\delta l.$$
Now, by direct computation:
$$\begin{array}{rl}
\lambda-1 &= WW'+W'W \vspace{2pt} \\
&=\alpha((X-tY)(\gamma X+\delta Y)+(\gamma X+\delta Y)(X-tY)) \vspace{2pt} \\
&=\alpha\gamma(2a-ts)+\alpha\delta(s-2t(l-b))\\
&=-t\gamma-\delta.
\end{array}$$
Let us first assume $\lambda=1$. Then, $\gamma(s-tl)=0$. If $\gamma=0,$ then $\delta=0$ and so $W'=0$. This means that
the $\phi(h)$-action is identically zero, yielding $s=l=0$. Otherwise, $s=tl$ and we are done. \par \smallskip

We finally show that the possibility $\lambda=-1$ can not occur. If $\lambda=-1$, then $\delta=2-t\gamma$ and $s\gamma=-(2-t\gamma)l$.
On the other hand,
\begin{equation}\label{sgamma}
l=\phi(h)\cdot (Y\# 1)=W'(Y\# 1)+(Y\# 1)W'=s\gamma +2(2-t\gamma)(l-b)
\end{equation}
Moreover,
$$\begin{array}{rl}
0 & = (W')^2 \vspace{2pt} \\
   & = \gamma^2 a +\delta^2(l-b) +\gamma \delta s \vspace{2pt} \\
   & \hspace{-10pt} \stackrel{(\ref{a=})}{=} \gamma^2(ts-t^2(l-b))+(2-t\gamma)^2(l-b)+\gamma(2-t\gamma)s \vspace{2pt} \\
   & = 2(l-b)(2-2t\gamma)+2\gamma s
\end{array}$$
From here, $s\gamma = (2t\gamma-2)(l-b).$ Substituting this in (\ref{sgamma}) we get $l=2b$, contradicting the fact that $C(b; 1,l)$ is $(H_4,R_l)$-Azumaya. \medskip

(2) If $l\not=0$, then $Im(\iota_l)=Im(i_{l^{-1}})$ by Proposition \ref{BM=BC} and the
statement follows from (1). It remains to show that $[C(a;t,s)]\in Im(\iota_0)$ implies $t=0$. If
$[C(a;t,s)] \in Im(\iota_0),$ there exists $b\in k^\cdot$ and an $H_4$-Azumaya algebra $A$
with trivial $h$-action and trivial $\phi(h)$-action such that $[C(a;t,s)]=[A\# C(b;0,1)].$ Then
$C(a;t,s)\# C(-b;0,1)\# \overline{A}\cong{\rm End}(P)$ for some
$D(H_4)$-module $P$. Arguing as in (1) we see that there is
$W=\alpha X\#1+\beta Y\# 1\in (C(a;t,s)\# C(-b;0,1))\# \overline{A}$
for some $\alpha,\beta\in k$ such that
$$\begin{array}{ll}
   & \hspace{11pt} h\cdot Z=W(g\cdot Z)-ZW,\\
0 & =h\cdot W=-2W^2=\alpha t+\beta,\\
t  & =h\cdot (X\# 1)=-2a\alpha,\\
0  & =h\cdot (Y\# 1)=2b\beta.
\end{array}
$$
From here if follows that $t=0$. \hfill$\Box$

\begin{corollary}\label{equal}
Let $[C(a;t,s)]$, $[C(b;p,q)]$ be in $BQ(k,H_4)$. Then $[C(a;t,s)]=[C(b;p,q)]$ if and only if
$C(a;t,s)\cong C(b;p,q)$.
\end{corollary}

\noindent \pf We analyze the case $t\neq0$, the other cases are treated
similarly. If $[C(a;t,s)]=[C(b;p,q)]$ and  $p=0$ then $[C(a;t,s)]\in
Im(\iota_0),$ contradicting Theorem \ref{distinct}. Then $tp\neq0$ and we
may reduce to the case $[C(a;1,s)]=[C(b;1,q)]\in Im(i_q)$. Applying again Theorem
\ref{distinct} we see that $s=q$ and the equality of classes is an
equality in $BM(k, H_4, R_q)$. Applying
$\Phi_0^{-1}\Psi_q^{-1}\Phi_q$ we obtain the equality
$[C(a-2^{-1}q;1,0)]=[C(b-2^{-1}q;1,0)]$ in $BM(k, H_4, R_0)$. From Proposition \ref{BM0},
we obtain $(4a-2q)^{-1}=(4b-2q)^{-1}$ and we have the statement.\hfill$\Box$

\begin{theorem}\label{intersection}
Let $i_t:BM(k, H_4, R_t)\to BQ(k, H_4)$ and $\iota_s:BC(k, H_4, r_s)\to BQ(k, H_4)$ be the natural
embeddings in $BQ(k, H_4)$. Then:
\begin{enumerate}
\item[(1)] $Im(i_t)\cap Im(\iota_s)\neq i_0(BW(k))$ if and only if $ts=1$. If this is the case, then $Im(i_t)=Im(\iota_s)$;
\item[(2)] $Im(i_t)\cap Im(i_s)\neq i_0(BW(k))$ if and only if $t=s$;
\item[(3)] $Im(\iota_t)\cap Im(\iota_s)\neq i_0(BW(k))$ if and only if $t=s$.
\end{enumerate}
\end{theorem}

\noindent \pf This is a consequence of Propositions \ref{BM0}, \ref{BCs}, \ref{BMt}, \ref{BM=BC} and Theorem \ref{distinct}. \hfill$\Box$

\section{The action of $Aut(H_4)$ on $Im(i_t)$ and $Im(\iota_s)$}

For a Hopf algebra $H$, a group morphism from ${\rm Aut}_{\rm Hopf}(H)$ to $BQ(k,H_4)$
has been constructed in \cite{CVZ2}, where the case of $H_4$ was also analized. The image of an
automorphism $\alpha$ can be represented as follows.

Let us denote by $H_\alpha$ the right $H$-comodule $H$ with left $H$-action
$l\cdot m=\alpha(\l2)mS^{-1}(\l1)$. Then $A_\alpha={\rm End}(H_\alpha)$ can be endowed of the $H$-Azumaya
algebra structure:
$$\begin{array}{l}
(l\cdot f)(m)=\l1\cdot f(S(\l2)\cdot m), \vspace{2pt} \\
\rho(f)(m)=\sum f(\m0)_{(0)}\otimes S^{-1}(\m1)f(\m0)_{(1)}.
\end{array}$$
The assignment $\alpha \mapsto [A_{\alpha^{-1}}]$ defines a
group morphism ${\rm Aut}_{\rm Hopf}(H)\to BQ(k, H)$.
The image of ${\rm Aut}_{\rm Hopf}(H)$ acts on $BQ(k, H)$ by conjugation. An easy
description of $[B(\alpha)]:=[A_\alpha][B][A_{\alpha}]^{-1}$ for any
representative $B$ has been given in \cite[Theorem 4.11]{CVZ2}. As an
algebra $B(\alpha)$ coincides with $B$, while the $H$-action and
$H$-coaction are:
\begin{equation}\label{autom}
h\cdot_{\alpha}b=\alpha(h)\cdot b,\quad\rho_{\alpha}(b)=\b0\otimes\alpha^{-1}(\b1).
\end{equation}

When $H=H_4$ the Hopf automorphism group is ${\rm Aut}_{\rm Hopf}(H_4)\cong k^\cdot$ and consists of the morphisms
that are the identity on $g$ and multiply $h$ by a nonzero scalar $\alpha$. The module $H_\alpha$ has action
$$\begin{array}{l}
g\cdot g=g, \quad g\cdot h=-h, \vspace{2pt} \\
h\cdot g=\alpha h g+g^2S^{-1}(h)=-(1+\alpha)gh, \quad h\cdot h=0,
\end{array}$$
and the kernel of the group morphism consists of $\{\pm 1\}$. We may thus embed
$(k^\cdot)^2\cong k^\cdot/\{\pm1\}$ into $BQ(k,   H_4)$ (cf. \cite{VOZ2}). We shall denote by $K$ the image of this group morphism. \medskip

We analyze this action on the classes and subgroups described in the
previous sections.

\begin{lemma}\label{conj-autom}
Let $\alpha\in k^\cdot$. Then:
\begin{enumerate}
\item[(1)] $[A_\alpha][C(a;t,s)][A_{\alpha}]^{-1}=[C(a;\alpha t,s\alpha^{-1})]$.
\item[(2)] $K$ acts trivially on $i_0(BW(k))$.
\end{enumerate}
In particular, $BM(k,H_4, R_{l\alpha^2})$ is conjugate to $BM(k,H_4, R_l)$ in $BQ(k, H_4)$ while $BM(k,H_4,R_0)$ and
$BC(k,H_4, r_0)$ are normalized by $K$.
\end{lemma}

\noindent \pf (1) It follows from direct computation that
$$h\cdot_{\alpha} x=\alpha t, \quad g\cdot_{\alpha}x=-x, \quad \rho(x)=x\otimes g+s\alpha^{-1}\otimes h.$$

(2) Since: the action of an automorphism of $H_4$ is trivial on $g$;
the action of $h$ is trivial on a representative of a class in
$BW(k)$; and the comodule map on a representative $A$ of a class in $BW(k)$ has image
in $A\otimes k{\mathbb Z}_2$, the formulas in (\ref{autom}) do not
modify the action and coaction on $A$ therefore $[A]=[A_\alpha][A][A_{\alpha}]^{-1}$ for every $[A]\in i_0(BW(k))$.
\smallskip

Since  $Im(i_l)$ is generated by $i_0(BW(k))$ and the classes $[C(a;1,l)]$, we see that $Im(i_l)$ is conjugate to
$Im(i_{\alpha^2l})$ in $BQ(k, H_4)$. If $l=0$ we get the statement concerning $Im(i_0)$. The statement concerning
$BC(k, H_4, r_0)$ follows because this group is generated by $i_0(BW(k))$ and the classes $[C(a;0,1)]$. \hfill$\Box$

\begin{remark}{\rm The observation that $Im(i_0)$ is normalized by $K$ has already been proved in \cite[\S 4]{VOZ4}. Lemma \ref{conj-autom} should be seen as a generalization of that result.}
\end{remark}

It is shown in \cite{Rad} that $(H_4,R_t)$ is equivalent to $(H_4,R_s)$ if and only if $t=\alpha^2s$ for some $\alpha\in k^\cdot$.
The above lemma shows that the Brauer groups of type $BM$ are conjugate in $BQ(k,H_4)$ if the corresponding triangular structures are equivalent. This is a general fact:

\begin{proposition}
Let $R$ and $R'$ be two equivalent quasitriangular structures on $H$ and let $\alpha\in{\rm Aut}_{\rm Hopf}(H)$ be such that $(\alpha\otimes\alpha)(R')=R$. Then the images of $BM(k,H,R)$ and $BM(k,H,R')$ are conjugate by the image of $\alpha$ in $BQ(k,H)$.
\end{proposition}

\noindent \pf If $B$ represents an element in $BM(k,H,R)$ then there will be an
action $\cdot$ on $B$ such that the coaction $\rho$ is given by $\rho(b)=(R^{(2)}\cdot b)\otimes
R^{(1)}$ for all $b \in B$. The image of $\alpha$ in $BQ(k,H)$ is represented by $A_{\alpha^{-1}}$. A representative of
$[A_{\alpha}]^{-1}[B][A_{\alpha}]$ is given by  the algebra $B$ with action
$h\cdot_{\alpha^{-1}}b=\alpha^{-1}(h)\cdot b$. The coaction is given by
$$\rho_\alpha(b)=(R^{(2)}\cdot b)\otimes\alpha(R^{(1)})=(\alpha(R^{(2)})\cdot_{\alpha}b)\otimes
\alpha(R^{(1)})=R'^{(2)}\cdot_{\alpha}b\otimes R'^{(1)},$$
so the coaction on $[A_{\alpha}]^{-1}[B][A_{\alpha}]$ is induced by
$R'$ and $\cdot_{\alpha}$.\hfill$\Box$ \medskip

For the dual statement, the proof is left to the reader.

\begin{proposition}Let $r$ and $r'$ be two equivalent coquasitriangular
  structures on $H$ and let $\alpha\in{\rm Aut}_{\rm Hopf}(H)$ be such
  that $r'(\alpha\otimes\alpha)=r$. Then the images of $BC(k,H,r)$
  and $BM(k,H,r')$ are conjugate by the image of $\alpha$ in $BQ(k,H)$.
\end{proposition}

\section{The subgroup $BQ_{grad}(k, H_4)$}

In this section we shall analyze the classes  that can be represented by $H_4$-Azumaya algebras for which the gradings coming from the $g$-action and the comodule structure coincide. They form a subgroup that will be related to the Brauer group $BM(k,E(2), R_N)$ of Nichols $8$-dimensional Hopf algebra $E(2)$ with respect to the  quasitriangular structure $R_N$ attached to the $2\times 2$-matrix $N$ with $1$ in the $(1,2)$-entry and zero elsewhere. \medskip

Let $BQ_{grad}(k,H_4)$ be the set of classes that can be represented by a $H_4$-Azumaya algebra $A$
for which the $|\cdot |$-grading and the $\deg$-grading coincide. In other words, the classes in $BQ_{grad}(k,H_4)$ can be represented by $D(H_4)$-module algebras on which the actions of $g$ and $\phi(g)$ coincide.  The last defining relation of $D(H_4)$ in Section 1 implies that the action of $h$ and $\phi(h)$ on such representatives commute. Clearly, $BQ_{grad}(k,H_4)$ is a subgroup of $BQ(k,H_4)$.

\begin{proposition}\label{subgroup}
$BQ_{grad}(k,H_4)$ is normalized by $K$.
\end{proposition}

\noindent \pf Let $[A]\in BQ_{grad}(k,H_4)$ with $|a|=\deg(a)$ for every $a\in A$ and let $[A_{\alpha}] \in K$. Then
$[A_\alpha\# A\# \overline{A_{\alpha}}]$ is represented by $A$ with action and coaction determined by (\ref{autom}). Since
$g$ is fixed by all Hopf automorphisms of $H_4$ we have
$$g\cdot_{\alpha}a=g\cdot a, \quad(\id\otimes\pi)\rho_\alpha(a)=(\id\otimes\pi)\rho(a),$$
so the two gradings are not modified by conjugation by $[A_\alpha]$. \hfill$\Box$ \medskip

The subgroup $BQ_{grad}(k,H_4)$ consists of those classes that can be represented by module algebras for the quotient of
$D(H_4)$ by the Hopf ideal $I$ generated by $\phi(g)\bowtie 1-\varepsilon\bowtie g$. Let us denote by $\pi_I$ the canonical
projection onto $D(H_4)/I$.

Let $E(2)$ be the Hopf algebra with generators $c,\,x_1,\,x_2,$ with relations
$$c^2=1,\quad x_i^2=0,\quad cx_i+x_ic=0,  \  i=1,2, \quad x_1x_2+x_2x_1=0,$$
coproduct
$$\Delta(c)=c\otimes c,\quad \Delta(x_i)=1\otimes x_i+x_i\otimes c,$$
and antipode $$S(c)=c,\quad S(x_i)=cx_i.$$

The Hopf algebra morphism
$$
\begin{array}{rl}
T\colon D(H_4)& \hspace{-5pt} \longrightarrow  E(2)\\
\phi(g)\bowtie 1& \hspace{-5pt} \mapsto c\\
\varepsilon\bowtie g & \hspace{-5pt} \mapsto c\\
\varepsilon\bowtie h & \hspace{-5pt} \mapsto x_1\\
\phi(h)\bowtie 1& \hspace{-5pt} \mapsto cx_2
\end{array}
$$
determines a Hopf algebra isomorphism $D(H_4)/I\cong E(2)$. The canonical quasitriangular structure ${\cal R}$ on $D(H_4)$ is
$$\begin{array}{rl}
{\cal R}&= \frac{1}{2}[\varepsilon\bowtie(1\otimes 1^*+g\otimes g^*+h\otimes h^*+gh\otimes (gh)^*)\bowtie 1] \vspace{3pt} \\
& \hspace{3pt} +\frac{1}{2}[\varepsilon\bowtie(1\otimes \varepsilon +g\otimes \varepsilon + 1\otimes \phi(g)-g\otimes\phi(g) \vspace{3pt} \\
& \hspace{20pt} +h\otimes \phi(h)+h\otimes\phi(gh)+gh\otimes\phi(h)-gh\otimes\phi(gh))\bowtie 1]
\end{array}$$
so $(\pi_I\otimes\pi_I)({\cal R})$ is a quasitriangular structure for $D(H_4)/I \cong E(2)$. Applying $T\otimes T$ to ${\cal R}$ we have:
\begin{equation}\label{RN}
\begin{array}{rl}
(T\otimes T)({\cal R}) &=\frac{1}{2}(1\otimes 1+1\otimes c+c\otimes
  1-c\otimes c\\
& \hspace{20pt} +x_1\otimes cx_2+x_1\otimes x_2+cx_1\otimes cx_2-cx_1\otimes x_2)
\end{array}
\end{equation}
The quasitriangular structures on $E(n)$ were computed in \cite{PVO1}. They are in bijection with $n\times n$-matrices with
entries in $k$. For a given matrix $M$ the corresponding quasitriangular structure is denoted by $R_M$. The map $T$ induces a quasitriangular morphism from $(D(H_4),{\cal R})$ onto $(E(2), R_N),$ where $N$ is the $2\times 2$-matrix with $1$ in the $(1,2)$-entry and zero elsewhere. If $A$ is a representative of a class in $BQ_{grad}(k,H_4)$ on which the ideal $I$ acts
trivially, then $A$ is an $E(2)$-module algebra and the maps $F$ and $G$ on $A\otimes A$ are the same as those induced by
$R_N$, so $A$ is $(E(2),R_N)$-Azumaya.

\begin{theorem}\label{sequence}The group $BM(k,E(2),R_N)$ fits into the following
exact sequence
$$\begin{array}{l}
\CD 1\longrightarrow{\mathbb Z}_2@>>>BM(k,
E(2),R_N)@>{T^*}>>BQ_{grad}(k,H_4)\longrightarrow 1.
\endCD
\end{array}$$
\end{theorem}

\noindent \pf Restriction of scalars through $T$ provides a group morphism $T^*$ from $BM(k,E(2), R_N)$ to $BQ(k,H)$ whose image is $BQ_{grad}(k,H_4)$. The kernel of $T^*$ consists of those classes $[A]$ such that $A\cong {\rm End}(P)$ as $D(H_4)$-module algebras, for some $D(H_4)$-module $P$. The class $[A]$ may be non-trivial only if $g$ and $\phi(g)$
act differently on $P$ even though they act equally on ${\rm End}(P)$. The $\phi(g)$- and $g$-action on ${\rm End}(P)$ are strongly
inner, hence there are elements $U$ and $u$ in ${\rm End}(P)$ such that $\phi(g)\cdot f=UfU^{-1}=ufu^{-1}=g\cdot f$
for every $f\in{\rm End}(P).$ Since ${\rm End}(P)$ is a central algebra, $U^2=u^2=1$, $uU=Uu$. From here, $U=\pm u,$ and if $[{\rm End}(P)]\neq 1$ in $BM(k,E(2),R_N)$ we necessarily have $U=-u$. The actions of $g$ and $\phi(g)$ on $P$ are given by the element $u$ and $U$ respectively, so for every non-trivial $[A]$ in ${\rm Ker}(T^*)$ we have $A\cong {\rm End}(P)$ for some $D(H_4)$-module $P$ for which $g$ acts as $-\phi(g)$. We claim that there is at most one non-trivial element in ${\rm Ker}(T^*)$.

Given any pair of such elements ${\rm End}(P)$ and ${\rm End}(Q)$ representing classes in ${\rm Ker}(T^*)$ we have ${\rm End}(P)\#{\rm End}(Q)\cong{\rm End}(P\otimes Q)$ as $D(H_4)$-module algebras by \cite[Proposition 4.3]{CVZ}, where $P\otimes Q$ is a $D(H_4)$-module. Then, the actions of $g$ and $\phi(g)$ on $P\otimes Q$ coincide, so $P\otimes Q$ is an $E(2)$-module. Thus, $[{\rm End}(P)][{\rm End}(Q)]$ is trivial in $BM(k,E(2), R_N)$ for every choice of $P$ and $Q$. Therefore, ${\rm Ker}(T^*)$ is either trivial or isomorphic to ${\mathbb Z}_2$. The proof is completed once we provide a non-trivial element. Let us consider $P=k^2$ on which $g,\,h,\,\phi(g)$ and $\phi(h)$ act via the following matrices $u,\,w,\,U,\,W$, respectively:
$$u=\left(\begin{array}{cc}
1&0\\
0&-1
\end{array}\right),\quad w=\left(\begin{array}{cc}
0&0\\
-2&0
\end{array}\right),\quad
U=-u,\,\quad W=\left(\begin{array}{cc}
0&1\\
0&0
\end{array}\right).$$
Then $P$ is a $D(H_4)$-module but not an $E(2)$-module. On the other hand, the $D(H_4)$-module algebra structure on
${\rm  End}(P)$ is in fact an $E(2)$-module algebra structure:
\begin{equation}\label{action1}
g\cdot f=u f u^{-1}=Uf U^{-1}=\phi(g)\cdot f;
\end{equation}
\begin{equation}\label{action2}
h\cdot f=w f u^{-1}+f uw, \quad
\phi(h)\cdot f=W f-U f U^{-1}W.
\end{equation}
Moreover, ${\rm  End}(P)$ is $(E(2),R_N)$-Azumaya because it is $H_4$-Azumaya.
We claim that the class of ${\rm  End}(P)$ is not trivial in
$BM(k,E(2),R_N)$. Indeed, if it were trivial, then the $E(2)$-action on
${\rm End}(P)$ given by $c.f=g.f$, $x_1.f=h.f$ and
$(cx_2).f=\phi(h).f$ would be strongly inner.
In other words, there would exist a convolution invertible algebra morphism $p\colon E(2)\to
{\rm End}(P)$ for which $l\cdot f=\sum p(\l1)f p^{-1}(\l2)$ for every
$l\in E(2)$. Putting $u'=p(c)$ we have
$c.f=u'f(u')^{-1}=ufu^{-1}$. Since  ${\rm End}(P)$ is a central
simple algebra, we necessarily have $u'=\lambda u$ and since $(u')^2=1$ we have $\lambda=\pm 1$.
Putting $w'=p(x_1)$ we have $x_1.f=w'fu'-fw'u'$ and since
$u'w'=-w'u'$,
we have $\lambda w'fu+\lambda fuw'=x_1.f=h.f=wfu+fuw$ for every $f\in
{\rm End}(P)$.  Using $uw=-wu$ we see that $(\lambda w'-w)f=f(\lambda
w'-w)$ so $w=\lambda w'+\mu$ for some $\mu\in k$. Using once more skew-commutativity of $u$ with $w$ and $w'$ we see that $\mu=0$.

Putting $W'=p(cx_2)$ and using that $u'W'=-W'u'$ we see that
$W'f-ufuW'=(cx_2).f=\phi(h).f=Wf-ufuW$ for every $f\in{\rm
  End}(P)$. From here, we deduce that $u(W'-W)=\nu\in k$. Using
skew-commutativity of  $u$ with $W$ and $W'$ we conclude that
$\nu=0$ so $W'=W$. Then $W'w'-w'W'=\lambda(Ww-wW)\neq 0$ so that relation
$(cx_2)x_1-x_1(cx_2)=0$ in $E(2)$ cannot be respected. Hence, $[{\rm End}(P)]\neq 1$ in $BM(k,E(2),R_N)$ and
${\rm Ker}(T^*)\cong{\mathbb Z}_2$.\hfill$\Box$ \medskip

The following proposition shows that the groups $BM(k,H_4,R_l)$ may be viewed inside $BM(k,E(2),R_N)$ and it also describes the image through $T^*$ of them.

\begin{proposition}\label{lambdamu}
For every $(\lambda,\mu)\in k\times k$ there is a group homomorphism $$\Theta_{\lambda,\mu}\colon
BM(k,H_4, R_{\lambda\mu})\to BM(k, E(2),R_N)$$
satisfying:
\begin{enumerate}
\item[(1)] The image of $\Theta_{0,0}$ is the subgroup isomorphic to $BW(k)$ represented by elements with
  trivial $x_1$- and $x_2$-action and $Ker(\Theta_{0,0}) \cong (k,+)$.

\item[(2)] $\Theta_{\lambda,\mu}$ is injective if and only if $(\lambda,\mu)\neq(0,0)$.

\item[(3)] For $(\lambda,\mu)\neq (0,0),$ the image of $T^*\Theta_{\lambda,\mu}$ is
$Im(i_{\mu\lambda^{-1}})$ if $\lambda\neq0$ and $Im(\iota_{\mu^{-1}\lambda})$ if $\mu\neq0$.
\end{enumerate}
\end{proposition}

\noindent \pf For every $(\lambda,\mu)\in k\times k$ the map $\theta_{\lambda,\mu}\colon E(2)\to H_4$ mapping $c\to g$,
$x_1\to \lambda h$ and $x_2\to \mu h$ is a Hopf algebra projection. A direct computation shows that
$(\theta_{\lambda,\mu}\otimes \theta_{\lambda,\mu})(R_N)=R_{\lambda\mu}$ so the pull-back of $\theta_{\lambda,\mu}$ induces the desired homomorphism $\Theta_{\lambda,\mu}$. \smallskip

(1) Let $(\lambda,\mu)=(0,0)$. Then any element in $BM(k, H_4, R_0)$ can be written as a pair of the form $([C(a;t,0)],[B])$ for $[B]\in BW(k)$. The image through $\Theta_{0,0}$ of such an element is $[C(a)][B]\in BW(k)$ with trivial $x_i$-action on $C(a)$. Clearly, $BW(k)=Im(\Theta_{0,0})$. That $Ker(\Theta_{0,0})$ is isomorphic to $(k,+)$ follows from the isomorphism
$BM(k, H_4, R_0)\cong(k,+)\times BW(k)$ and the fact that $(k,+)$ is realized as classes admitting a representative that is trivial when viewed as a $k{\mathbb Z}_2$-module algebra. \smallskip

(2) Let $(\lambda,\mu)\neq(0,0)$. If $\Theta_{\lambda,\mu}([A])=1$
then $A$ is isomorphic to an endomorphism algebra with strongly inner
$E(2)$-action. In other words, $A\cong {\rm End}(P)$ and there is a convolution invertible algebra map $p\colon E(2)\to A$ such that  $l\cdot a =\sum p(\l1)ap^{-1}(\l2)$ for every $l\in E(2), a\in A$. There are elements $u,v,w\in A$
with $u$ invertible such that $c\cdot a=g\cdot a= uau^{-1}$, $x_1\cdot a=(v a-av)u=\lambda h\cdot a$ and  $x_2\cdot a=(wa-aw)u=\mu h\cdot a$. Then $$0=\mu x_1\cdot a-\lambda x_2\cdot a=((\mu v-\lambda w)a-a(\mu v-\lambda w))u \quad \forall a\in A,$$ and since $u$ is invertible and $A$ is central we have $\mu v-\lambda w=\eta$ for some $\eta \in k$. The relation between $v$ and $w$ gives $\eta=0$ and so $\mu v=\lambda w.$ Thus, the same elements $u,v$ and $w$ ensure that the $H_4$-action on $A$ is strongly inner. Therefore $[A]=1$ in $BM(k,H_4, R_{\lambda\mu})$. The converse follows from (1). \smallskip

(3) Let us now assume that $(\lambda,\mu)\neq(0,0).$ It is immediate to see that if $[A]\in BW(k)\subset
BM(k,H_4, R_{\lambda\mu})$ is represented by an algebra with trivial $h$-action, then $\Theta_{\lambda,\mu}([A])$ is
represented by an algebra with trivial $x_1$- and $x_2$-action. Hence $T^*\Theta_{\lambda,\mu}(BM(k, H_4, R_{\lambda\mu}))\subset i_0(BW(k))$ and the restriction of $T^*\Theta_{\lambda,\mu}$ to $BW(k)$ is an isomorphism onto
$i_0(BW(k))$. Let us now consider the class $[C(a;1,\lambda\mu)]\in BM(k, H_4, R_{\lambda\mu})$. Its image
through $\Theta_{\lambda,\mu}$ is the algebra generated by $x$ with
$x^2=a$, with $c\cdot x=-x$, $x_1\cdot x=\lambda$ and $x_2\cdot
x=\mu$. A direct verification shows that
$T^*\Theta_{\lambda,\mu}([C(a;1,\lambda\mu)])=[C(a;\lambda,\mu)]$.
Then the image of $T^*\Theta_{\lambda,\mu}$ is
$Im(i_{\mu\lambda^{-1}})$ if $\lambda\neq0$ and
$Im(\iota_{\mu^{-1}\lambda})$ if $\mu\neq0$. \hfill$\Box$ \bigskip

Theorem \ref{sequence} shows that one should  understand $BM(k,E(2),R_N)$ in order to compute $BQ(k,H_4)$. In view of Proposition \ref{lambdamu}, $BM(k,E(2),R_N)$ seems to be much more complex that the groups of type BM treated in \cite{hnu,enne,VOZ3}.

\section{Appendix}

This last section is devoted to the analysis of some difficulties occurring in the study of the structure of  $(E(2),R_N)$-Azumaya algebras. We show that the set of classes represented by ${\mathbb Z}_2$-graded central simple algebras (with respect to the grading induced by the $c$-action) is not a subgroup of  $BM(k,E(2),R_N)$. \medskip

Let us consider the braiding $\psi_{VW}$ determined by $R_N$ between two left $E(2)$-modules $V$ and $W$. Let $v \in V$ and $w \in W$ be homogeneous elements with respect to the ${\mathbb Z}_2$-grading induced by the $c$-action. By direct computation it is:
$$
\begin{array}{l}
\psi_{VW}(v\otimes w) =\sum R_N^{(2)}\cdot w\otimes R_N^{(1)}\cdot v \vspace{2pt} \\
\hspace{1.7cm}  = (-1)^{|v||w|}w\otimes v+(-1)^{|w|+1}(-1)^{(|v|+1)(|w|+1)}(x_2\cdot w)\otimes (x_1\cdot v).
\end{array}
$$
If we denote by $\psi_0$ the braiding associated with the ${\mathbb Z}_2$-grading we have
\begin{equation}\label{psi}
\psi_{VW}(v\otimes w)=\psi_0(v\otimes w)+(-1)^{|w|+1}\psi_0(x_1\cdot
v\otimes x_2\cdot w).
\end{equation}
Let $F$ and $G$ be the maps in \eqref{f-and-g} defining an $(E(2), R_N)$-Azumaya algebra $A$
and let $F_0$ and $G_0$ be the maps defining an $(E(2),R_0)$-Azumaya
algebra, that is, the maps determining when an $E(2)$-module algebra
is ${\mathbb Z}_2$-graded central simple. It is not hard to verify by
direct computation that, for homogeneous $a,\,b,\,d\in A$ with respect
to the $c$-action we have:
\begin{equation}\label{effe}
F(a\# b)(d)=F_0(a\# b)(d)+(-1)^{|d|+1}F_0(a\# x_1\cdot b)(x_2\cdot d)
\end{equation}
\begin{equation}\label{gi}
G(a\# b)(d)=G_0(a\# b)(d)+(-1)^{|a|+1}F_0(x_2\cdot a\# b)(x_1\cdot d)
\end{equation}

Notice that if either $x_1$ or $x_2$ acts trivially, then $F=F_0$ and $G=G_0$. So in this case, $A$ is $(E(2),R_N)$-Azumaya if and only if it is ${\mathbb Z}_2$-graded central simple (i.e. $A$ is $(E(2),R_0)$-Azumaya). We will say that the $x_i$-action on an $E(2)$-module algebra $A$ is {\em inner} if there exists an odd element $v\in A$ such that $x_i\cdot a= v(c\cdot a)-av$ for every $a\in A$.

\begin{theorem}\label{conditions}
Let $A$ be an $(E(2),R_N)$-Azumaya algebra. The following assertions are equivalent:
\begin{enumerate}
\item[(1)] The $x_1$-action on $A$ is inner;
\item[(2)] The $x_2$-action on $A$ is inner;
\item[(3)] $A$ is a ${\mathbb Z}_2$-graded central simple algebra.
\end{enumerate}
In addition, the $E(2)$-action on $A$ is inner if and only if  $A$ is a central simple algebra.
\end{theorem}

\noindent \pf (1) $\Rightarrow$ (3) Let $v_1\in A$ be an odd element such that $x_1\cdot a=v_1(c\cdot a)-av_1$ for all $a \in A$. Applying equality (\ref{effe}) to any homogeneous $b$ and $d$ in $A$ gives:
\begin{equation}\label{effe-inner}
\begin{array}{ll}
F(a\# b)(d) & =F_0(a\# b)(d)+F_0(a\# b)((x_2\cdot d)v_1) \vspace{2pt} \\
 & \hspace{10pt} + (-1)^{|d|}F_0(a\# bv_1)(x_2\cdot d)
\end{array}
\end{equation}
This equality extends to all elements $a$ and $b$ in $A$. If $A$ were not ${\mathbb Z}_2$-graded central simple, there would exist an element $0 \neq \sum_i a_i\# b_i$  in $Ker(F_0)$. Then $(\sum_i a_i\# b_i)(1\# v_1)=\sum_i a_i\# b_iv_1 \in Ker(F_0)$ and for every $f$ in $A$ we would have $F_0(\sum_i a_i\#b_i)(f)=F_0(\sum_i a_i\# b_iv_1)(f)=0$. It follows from (\ref{effe-inner}) that $\sum_i a_i\# b_i  \in Ker(F),$ contradicting the injectivity of $F$. \smallskip

(2) $\Rightarrow$ (3) Similarly to (1) $\Rightarrow$ (3) replacing $F$ by $G$.  \smallskip

(3) $\Rightarrow$ (1), (2) Suppose that $A$ is a ${\mathbb Z}_2$-graded central simple
algebra. If $A$ is a central simple algebra then the $E(2)$-action on $A$ is inner by the Skolem-Noether theorem.
If $A$ is not central simple then it is of odd type (\cite[Theorem  3.4, Definition 3.5]{lam}) and it is $(H_4,R_0)$-Azumaya for the subalgebra of $E(2)$, isomorphic to $H_4$ generated by $c$ and $x_i$. By \cite[Theorem 3.4]{ACZ} the $x_i$-action is inner. \smallskip

Let us finally assume that the $E(2)$-action on $A$ is inner. Then $A$ is
a ${\mathbb Z}_2$-graded central simple algebra. Since $E(2)$ acts innerly on $A$ then it acts trivially on its center $Z(A)$. Besides it is immediately seen that $Z(A)$ is contained in the right and left
$E(2)$-center, that are trivial because $A$ is assumed to be $E(2)$-Azumaya. Hence $Z(A)$ must
be trivial and so $A$ is also a central algebra. By the structure theorems
of ${\mathbb Z}_2$-graded central simple algebras (\cite[Theorem IV.3.4]{lam}), $A$ is central
simple.\hfill$\Box$

\begin{proposition}\label{independent}
Let $A$ and $B$ be two equivalent $(E(2),R_N)$-Azumaya algebras. Then the $x_i$-action on $A$ is inner if and only if it
is so on $B$.
\end{proposition}

\noindent \pf Let $P$ and $Q$ be finite dimensional $E(2)$-modules for which $A\# {\rm End}(P)\cong
B\# {\rm End}(Q)$. If the $x_i$-action on $A$ is inner then it is so on $A\# {\rm End}(P)$ by
\cite[Proposition 4.6]{enne}, hence it is so on $B\# {\rm End}(Q)$, which is a ${\mathbb Z}_2$-graded central
simple algebra by Theorem \ref{conditions}. For $i=1,2$, let $W_i,v_i$ be odd elements in $B\# {\rm End}(Q)$ and ${\rm End}(Q)$ respectively inducing the $x_i$-action. We recall that $x_j\cdot v_i=0$ because the action on ${\rm End}(Q)$ is strongly
inner, while $x_j\cdot W_i$ is a scalar for every pair $i,j$ because $x_j\cdot W_i$ belongs to the graded center of $B\# {\rm End}(Q)$. The odd elements
$T_i=W_i-1\# v_i-(x_2\cdot W_i)(1\# v_1)\in B\# {\rm End}(Q)$ for
$i=1,2$ are such that $x_j\cdot T_i=x_j\cdot W_i$ for every $i$ and
$j$. Moreover, for every homogeneous $f\in {\rm End}(Q)$ with respect
to the $c$-action we have:
$$\begin{array}{l}
(-1)^{|f|}T_i(1\# f) =W_i(c\cdot 1\# c\cdot f)-1\# v_i(c\cdot f)-(x_2\cdot W_i)(1\# v_1 (c\cdot f)) \vspace{2pt} \\
\hspace{1.3cm} =(1\# f)W_i-(1\# fv_i)-(x_2\cdot W_i)(1\# fv_1)-(x_2\cdot W_i)(x_1\cdot(1\# f)) \vspace{2pt} \\
\hspace{1.3cm} =(1\# f)[W_i-1\# v_i-(x_2\cdot W_i)(1\# v_1)]- (x_2\cdot W_i)(x_1\cdot(1\# f)) \vspace{2pt} \\
\hspace{1.3cm} =(1\# f)T_i-(x_2\cdot W_i)(x_1\cdot(1\# f)).
\end{array}$$
In other words,
$$(1\# f)T_i=(-1)^{|f||T_i|}T_i(1\# f)+(x_2\cdot T_i)(x_1\cdot (1\#
f)),$$ so by (\ref{psi}) the element $T_i\in C^l_{B\# {\rm End}(Q)}({\rm End}(Q)),$ the left centralizer
of ${\rm End}(Q)$ in $B\# {\rm End}(Q)$, that is, $T_i\in B\#1$ by the double centralizer theorem \cite[Theorem 2.3]{ACZ}.
Besides, for every homogeneous $b\in B$ we have:
$$\begin{array}{rl}
T_i(c\cdot b\# 1)-(b\# 1)T_i & = (-1)^{|b|}W_i(b\#1)-(b\# v_i)-(x_2\cdot W_i)(b\# v_1) \vspace{2pt} \\
 &  \hspace{10pt} -(b\#1)W_i+(b\# v_i)+(x_2\cdot W_i)(b\# v_1) \vspace{2pt} \\
& =x_i\cdot(b\# 1).
\end{array}$$
Hence the $x_i$-action on $B$ is inner.\hfill$\Box$ \medskip

We conclude by showing that, contrarily to the cases treated in the literature (\cite{hnu,enne,VOZ3}), a
Skolem-Noether-like approach  is probably
not appropriate for the computation of $BM(k,E(2),R_N)$ because the set of classes admitting a representative with inner action is not
a subgroup.

\begin{theorem}\label{not-subgroup}
The classes in $BM(k,E(2), R_N)$ that are represented by ${\mathbb Z}_2$-graded central simple algebras
do not form a subgroup.
\end{theorem}

\pf Let $t\neq 0, 1$ and $q\neq 2$ be in $k$. We consider the representative $C(1;t,2)$ generated by $x$ with $x^2=1$, $c\cdot x=-x$, $x_1\cdot x=t$ and $x_2\cdot x=2$ and the representative $C(1;1,q)$ generated by $y$ with $y^2=1$, $c\cdot y=-y$, $x_1\cdot
y=1$ and $x_2\cdot y=q$. Both are $(E(2),R_N)$-Azumaya because $C(1;1,2t)$ is $(H_4,R_{2t})$-Azumaya, $C(1;1,q)$ is $(H_4,R_{q})$-Azumaya and $C(1;t,2), C(1;1,q)$ are obtained from these ones respectively by pulling back through $\theta_{\lambda,\mu}$. They are also ${\mathbb Z}_2$-graded central simple algebras. Their product
$C(1;t,2)\# C(1;1,q)$ is generated by the odd elements $X$ and $Y$ with $X^2=1$, $Y^2=1$ and $XY+YX=2$. The element $X-Y$ is easily seen to lie in the ${\mathbb Z}_2$-graded center, so $C(1;t,2)\# C(1;1,q)$ is not a ${\mathbb Z}_2$-graded central simple algebra. If $B$ were another representative of $[C(1;t,2)\# C(1;1,q)]$ that is a ${\mathbb Z}_2$-graded central simple algebra,
then by Theorem \ref{conditions}, the $x_1$-action on it would be inner. By Proposition \ref{independent}, $x_1$ would act innerly on $C(1;t,2)\# C(1;1,q)$. Applying again Theorem \ref{conditions}, $C(1;t,2)\# C(1;1,q)$ would be ${\mathbb Z}_2$-graded central simple. \hfill$\Box$ \vspace{1cm}

\begin{center}
{\bf Acknowledgements}
\end{center}

This research was partially supported by the Azioni Integrate Italia-Espa\~na AIIS05E34A {\it Algebre, coalgebre, algebre di Hopf e loro rappresentazioni}. The second named author is also supported by projects MTM2008-03339 from MCI and FEDER and P07-FQM-03128 from Junta de Andaluc\'{\i}a.

\end{document}